\documentclass[twoside]{article} 
\usepackage{graphicx}
\date{} 
\title{An asymptotic approximation for the Riemann zeta function revisited}
\author{\sc R. B.\ Paris \\
{\em Division of Computing and Mathematics,} \\
{\em Abertay University, Dundee DD1 1HG, UK}}
\begin{document}
\def\f#1#2{\mbox{${\textstyle \frac{#1}{#2}}$}}
\def\dfrac#1#2{\displaystyle{\frac{#1}{#2}}}
\def\boldal{\mbox{\boldmath $\alpha$}}
\newcommand{\bee}{\begin{equation}}
\newcommand{\ee}{\end{equation}}
\newcommand{\lam}{\lambda}
\newcommand{\ka}{\kappa}
\newcommand{\al}{\alpha}
\newcommand{\ba}{\beta}
\newcommand{\la}{\lambda}
\newcommand{\ga}{\gamma}
\newcommand{\eps}{\epsilon}
\newcommand{\om}{\omega}
\newcommand{\fr}{\frac{1}{2}}
\newcommand{\fs}{\f{1}{2}}
\newcommand{\g}{\Gamma}
\newcommand{\br}{\biggr}
\newcommand{\bl}{\biggl}
\newcommand{\ra}{\rightarrow}
\newcommand{\gtwid}{\raisebox{-.8ex}{\mbox{$\stackrel{\textstyle >}{\sim}$}}}
\newcommand{\ltwid}{\raisebox{-.8ex}{\mbox{$\stackrel{\textstyle <}{\sim}$}}}
\renewcommand{\topfraction}{0.9}
\renewcommand{\bottomfraction}{0.9}
\renewcommand{\textfraction}{0.05}
\newcommand{\mcol}{\multicolumn}
\date{}
\maketitle
\pagestyle{myheadings}
\markboth{\hfill \sc R. B.\ Paris  \hfill}
{\hfill \sc Asymptotic representation of $Z(t)$ \hfill}
\begin{abstract}
We revisit a representation for the Riemann zeta function $\zeta(s)$ expressed in terms of normalised incomplete gamma functions given by the author and S. Cang in Methods Appl. Anal. {\bf 4} (1997) 449--470. Use of the uniform asymptotics of the incomplete gamma function produces an asymptotic-like expansion for $\zeta(s)$ on the critical line $s=1/2+it$ as $t\to+\infty$. The main term involves the original Dirichlet series smoothed by a complementary error function of appropriate argument together with a series of correction terms. It is the aim here to present these correction terms in a more user-friendly format by expressing then in inverse powers of $\omega$, where $\omega^2=\pi s/(2i)$, multiplied by coefficients involving trigonometric functions of argument $\omega$.
\vspace{0.3cm}

\vspace{0.1cm}
 
\end{abstract}

\vspace{0.3cm}

\noindent $\,$\hrulefill $\,$

\vspace{0.3cm}

\begin{center}
{\bf 1.\ Introduction}
\end{center}
\setcounter{section}{1}
\setcounter{equation}{0}
\renewcommand{\theequation}{\arabic{section}.\arabic{equation}}
The computation of the Riemann zeta function on the critical line $s=\fs+it$, $t\geq 0$ is normally carried out using the real function $Z(t)=e^{i\vartheta(t)} \zeta(\fs+it)$, where the phase angle
\[\vartheta(t)=\arg\,\g(\f{1}{4}+\fs it)-\fs t \log\,\pi.\]
In \cite{PC}, a representation of $Z(t)$ on the critical line was given in the form
\bee\label{e11}
Z(t)=2 \Re e^{i\vartheta(t)} \bl\{\sum_{n=1}^\infty n^{-s} Q(\fs s, \pi n^2 i)-\frac{\pi^{s/2} e^{\pi is/4}}{s \g(\fs s)}\br\},
\ee
wherein $s=\fs+it$, $t\geq 0$ and
$Q(a,z)=\g(a,z)/\g(a)$ is the normalised (upper) incomplete gamma function.  The behaviour of $Q(\fs s,\pi n^2 i)$ for large $t$ changes abruptly in the neighbourhood of its transition point $\fs s=\pi n^2i$; that is, when the summation index $n$ roughly  equals the Riemann-Siegel cut-off value $N_t$ given by
\bee\label{e12}
N_t=[(t/2\pi)^{1/2}],
\ee
where square brackets denote the integer part. Then for large $t$, the function $|Q(\fs s,\pi n^2i)|\simeq 1$ when $n\,\ltwid\,N_t$ and decays (algebraically) to zero when $n\,\gtwid\,N_t$. As a consequence, the absolutely convergent sum\footnote{The sum in (\ref{e11}) converges absolutely since $\g(a,z)\sim z^{a-1}e^{-z}$ as $|z|\to\infty$ in $|\arg\,z|<3\pi/4$, and so late terms behave like $(-)^n/n^2$.} in (\ref{e11}) represents the original Dirichlet series ``smoothed'' by the incomplete gamma function.

The uniform asymptotic expansion of $Q(a,z)$ valid for $a\to\infty$ when $|z|\in[0,\infty)$ in the domains $|\arg\,z|<\pi$ and $|\arg\,(z/a)|<2\pi$ was employed in \cite{PC} to derive an asymptotic approximation for $Z(t)$ valid as $t\to\infty$. The resulting expansion involved the Dirichlet sum smoothed by a complementary error function of appropriate argument together with asymptotic correction terms in inverse powers of $\om$, where $\om^2=\pi s/(2i)$, decorated by functions $S_k(\om)$ involving $\csc \om$ and its derivatives. In \cite{P} an attempt was made to establish the asymptotic nature of this expansion; this was not successful, however, due to the use of an insufficiently precise error bound for the incomplete gamma function.

In \cite{PC} the quantities $S_k(\om)$ were treated as computable functions.
The purpose of the present paper is to revisit the expansion (\ref{e11}) with the aim of presenting the coefficients in the correction terms in a more user-friendly format. This is achieved using the expansion of the functions $S_k(\om)$ as polynomials in $1/\om$ of degree $k-1$. 
\vspace{0.6cm}

\begin{center}
{\bf 2.\ The asymptotic approximation for $Z(t)$}
\end{center}
\setcounter{section}{2}
\setcounter{equation}{0}
\renewcommand{\theequation}{\arabic{section}.\arabic{equation}}
To make this paper self-contained, we summarise in this section the main steps in the derivation of the asymptotic approximation derived from (\ref{e11}) and presented in \cite{PC}.
We employ the uniform asymptotic expansion of the normalised incomplete gamma function valid as $a\to\infty$ in $|\arg\,a|<\pi$ in the form \cite[p.~181]{DLMF}
\bee\label{e21}
Q(a,z)=\frac{1}{2} \mbox{erfc}\,(\eta\sqrt{a/2})+\frac{e^{-\frac{1}{2}a\eta^2}}{\sqrt{2\pi a}}\bl\{\sum_{k=0}^{m-1}c_k(\eta) a^{-k}+a^{-m}R_m(a,\eta)\br\}
\ee
for $m=1, 2, \ldots $, where
\bee\label{e22}
\fs\eta^2=\lambda-1-\log\,\lambda,\quad \lambda=\frac{z}{a},\quad \mu=\lambda-1
\ee
and the quantity $R_m(a,\eta)$ is the remainder term in the expansion truncated after $m$ terms.
The choice of the square root branch for $\eta(\lambda)$ is made such that $\eta(\lambda)$ and $\lambda-1$ have the same sign when $\lambda>0$; we then have $\eta\simeq\lambda-1$ when $\lambda\simeq 1$. The coefficients $c_k(\eta)$ are given by
\bee\label{e22a}
c_k(\eta)=\frac{(-)^kW_k(\mu)}{\mu^{2k+1}}-\frac{(-2)^k (\fs)_k}{\eta^{2k+1}},\qquad W_k(\mu):=\sum_{r=0}^{2k}\al_{r,k} \mu^r,
\ee 
where $W_k(\mu)$ is a polynomial in $\mu$ of degree $2k$ and the coefficients $\al_{r,k}$ are discussed in \cite[Appendix]{PC}; see also \cite{T}. The first few coefficients are given by
\begin{eqnarray*}
\al_{00}&=&1\\
\al_{01}&=& 1,\ \ \al_{11}=1,\ \ \al_{21}=\f{1}{12}\\
\al_{02}&=& 3,\ \ \al_{12}=5,\ \ \al_{22}=\f{25}{12},\ \ \al_{32}=\f{1}{12},\ \ \al_{42}=\f{1}{288}\\
\al_{03}&=&15,\ \al_{13}=35,\ \ \al_{23}=\f{105}{4},\ \ \al_{33}=\f{77}{12},\ \ \al_{43}=\f{49}{288},\ \ \al_{53}=\f{1}{288},\ \ \al_{63}=-\f{139}{51840}. 
\end{eqnarray*}

The modified complementary error function is introduced by
\bee\label{e23}
\mbox{erfc}\,(z;m)=\mbox{erfc}\,z-\frac{e^{-z^2}}{\sqrt{\pi}} \sum_{r=0}^{m-1}\frac{(-)^r (\fs)_r}{z^{2r+1}}\qquad (z\neq 0;\ m=1, 2, \ldots\,),
\ee
which corresponds to the deletion from $\mbox{erfc}\,z$ of the first $m$ terms of its asymptotic expansion for $|z|\to\infty$ in $|\arg\,z|<\f{3}{4}\pi$. Then the expansion (\ref{e21}) can rewritten in the form
\bee\label{e24}
Q(a,z)=\frac{1}{2} \mbox{erfc}\,(\eta\sqrt{a/2};m)+\frac{e^{-\frac{1}{2}a\eta^2}}{\sqrt{2\pi a}}\bl\{\sum_{k=0}^{m-1}\frac{(-)^kW_k(\mu)}{\mu^{2k+1}} a^{-k}+a^{-m}R_m(a,\eta)\br\}.
\ee

We now substitute the expansion (\ref{e24}) into (\ref{e11}) when we identify the parameters $a$, $\lambda\equiv\lambda_n$ and $\mu\equiv \mu_n$ (and correspondingly $\eta\equiv \eta_n$) with
\bee\label{e25}
a=\f{1}{2}s=\f{1}{4}+\f{1}{2}it,\quad \lambda_n=2\pi n^2 i/s,\quad \mu_n=(\pi n/\om)^2-1,\quad \om^2=\pi s/(2i).
\ee
We also make use of the well-known expansion 
\[\frac{1}{\g(\fs s)}=(2\pi)^{-1/2} (\fs s)^{(1-s)/2} e^{s/2} \bl\{\sum_{r=0}^{m-1} \gamma_r (\fs s)^{-r}+ (\fs s)^{-m} H_m(s)\br\}\qquad (t\to\infty),\]
where the $\gamma_r$ are the Stirling coefficients with $\gamma_1=1$, $\gamma_2=-\f{1}{12}$, $\gamma_2=\f{1}{288}$, $\gamma_3=\f{139}{51840}, \ldots\ ,$ and $H_m(s)$ is a remainder term.  Since
\[e^{-\frac{1}{2}a\eta_n^2}=(-)^n n^s e^{\pi i s/4} (2\pi e/s)^{s/2},\] 
we then find that
\bee\label{e26}
Z(t)=\Re e^{i\vartheta(t)} \bl\{\sum_{n=1}^\infty n^{-s} \mbox{erfc}\,(\fs \eta_n \sqrt{s};m)-\frac{e^{\pi is/4}}{\sqrt{\pi s}}\,\bl(\frac{2\pi e}{s}\br)^{s/2}\bl( \sum_{r=0}^{m-1} \frac{(-)^rA_r(s)}{(\fs s)^r}+(\fs s)^{-m} {\hat R}_m\br)\br\},
\ee
where ${\hat R}_m$ is the remainder that contains contributions from the remainder term in $Q(a,z)$ and that in the expansion of $1/\g(\fs s)$. We do not consider this term here; an attempt was made in \cite{P} to bound ${\hat R}_m$ and demonstrate that (\ref{e26}) is an asymptotic expansion, but the bounds employed on the incomplete gamma function were not sharp enough to achieve this.

The coefficients $A_r(s)$ are given by
\[A_r(s)=(-)^r\gamma_r-2\sum_{n=1}^\infty \frac{(-)^n W_r(\mu_n)}{\mu_n^{2r+1}}=\sum_{n=-\infty}^\infty \frac{(-)^{n-1} W_r(\mu_n)}{\mu_n^{2r+1}}\]
where $W_r(\mu_0)=W_r(-1)=(-)^r\gamma_r$; see \cite[(A.1)]{PC}. From the definitions of $W_k(\mu_n)$ and $\mu_n$ in (\ref{e22a}) and (\ref{e25}), we then find that
\[A_r(s)=\sum_{k=0}^{2r}(-)^k \al_{r,k} \sum_{n=-\infty}^\infty (-)^n \bl(\frac{\om^2}{\om^2-(\pi n)^2}\br)^{\!2r+1-k}.\]
The inner sum may be expressed in terms of the functions $S_k(\om)$ defined by
\bee\label{e26a}
S_k(\om):=2^{k-1} \sum_{n=-\infty}^\infty \frac{(-)^n \om^k}{(\om^2-(\pi n)^2)^k}\qquad k=1, 2, \ldots 
\ee
so that
\bee\label{e27}
A_r(s)=2^{-2r}\om^{2r+1} B_r(\om),\qquad B_r(\om)=\sum_{k=0}^{2r} (-)^k \al_{r,k} \,\frac{S_{2r+1-k}(\om)}{(\om/2)^k}.
\ee
After some routine algebra we finally obtain the expansion on the critical line for large $t$ in the form \cite{PC}:
\newtheorem{theorem}{Theorem}
\begin{theorem}$\!\!\!.$ \ \ As $t\to+\infty$ on the critical line $s=\fs+it$, we have the expansion 
\bee\label{e28}
Z(t)=\Re e^{i\vartheta(t)}\bl\{\sum_{n=1}^\infty n^{-s} \mbox{erfc}\,(\fs \eta_n \sqrt{s};m)
-\frac{1}{\sqrt{2i}}\bl(\frac{\pi e^{1/2}}{\om}\br)^{\!s}\bl(T_m(\om)+\om^{-2m+1} {\hat R}_m\br)\br\},
\ee
where 
$\om^2=\pi s/(2i)$ and
\bee\label{e29}
T_m(\om):=\sum_{r=0}^{m-1} (\pi i/4)^r B_r(\om)
\ee
with the coefficients $B_r(\om)$ defined in (\ref{e27}).
\end{theorem}
In \cite{PC} the functions $S_k(\om)$ appearing in the coefficients $B_r(\om)$ were treated as computable functions. In the next section we express $B_r(\om)$, and hence $T_m(\om)$, as a finite series in inverse powers of $\om$ with coefficients containing  certain trigonometric functions.

\vspace{0.6cm}

\begin{center}
{\bf 3.\ The coefficients in the expansion}
\end{center}
\setcounter{section}{3}
\setcounter{equation}{0}
\renewcommand{\theequation}{\arabic{section}.\arabic{equation}}
In this section we examine the correction term $T_m(\om)$ defined in (\ref{e29})
and express it as a finite series in descending powers of $\om$.  To achieve this
we first observe that the functions $S_k(\om)$ defined in (\ref{e26a}) satisfy 
\[S_1(\om)=\csc\,\om, \qquad S_{k+1}(\om)=\frac{1}{\om} S_k(\om)-\frac{1}{k} \,\frac{dS_k(\om)}{d\om}\quad (k\geq 1).\]
From this we find with the help of {\it Mathematica} that
\bee\label{e300}
S_{k+1}(\om)=\frac{\csc \om}{k!} \sum_{r=0}^k \frac{b_{r,k}}{\om^k},
\ee
where the coefficients $b_{r,k}$ for $0\leq r,k\leq 6$ are
\begin{eqnarray*}
b_{00}&=&1; \\
\\[-0.30cm]
b_{01}&=&\cot\,\om,\ \ b_{11}=1;\\
\\[-.30cm]
b_{02}&=&\cot^2 \om+\csc^2 \om,\ \ b_{12}=3\cot \om,\ \ b_{22}=3;\\
\\[-0.30cm]
b_{03}&=&\cot^3 \om+5\cot \om \csc^2 \om, \ \ b_{13}=6(\cot^2 \om+\csc^2 \om),\ \ b_{23}=15\cot \om,\ \ b_{33}=15;\\
\\[-0.30cm]
b_{04}&=&\cot^4 \om+18\cot^2 \om \csc^2 \om+5\csc^4 \om,\ \ b_{14}=10(\cot^3 \om+5\cot \om \csc^2 \om),\\
b_{24}&=&45(\cot^2 \om+\csc^2 \om),\ \ b_{34}=105\cot \om,\ \ b_{44}=105;\\
\\[-0.30cm]
b_{05}&=&\cot^5 \om+58\cot^3 \om \csc^2 \om+61\cot \om \csc^4 \om,\ \ b_{15}=15(\cot^4 \om+18\cot^2 \om \csc^2 \om+5\csc^4 \om), \\
b_{25}&=&105(\cot^3 \om+5\cot \om \csc^2 \om),\ \ b_{35}=420(\cot^2 \om+\csc^2 \om),\ \ b_{45}=945\cot \om,\ \\\ b_{55}&=&945;\\
\\[-0.30cm]
b_{06}&=&\cot^6 \om+179\cot^4 \om \csc^2 \om +479\cot^2 \om \csc^4 \om+61\csc^6 \om,\\ 
b_{16}&=&21(\cot^5 \om+58\cot^3 \om \csc^2 \om+61\cot \om \csc^4 \om),\ \\
b_{26}&=&210(\cot^4 \om+18\cot^2 \om \csc^2 \om+5\csc^4 \om),\ \ b_{36}=1260(\cot^3 \om+5\cot \om\csc^2 \om),\  \\ b_{46}&=&4725(\cot^2 \om+\csc^2 \om), \ \ b_{56}=10395\cot \om,\ \ b_{66}=10395.\\
\end{eqnarray*}

Then, from (\ref{e27}) and (\ref{e300}), we find $B_0(\om)=\csc \om$ and
\[B_1(\om)=\al_{10}S_3(\om)-\frac{2}{\om}\al_{11} S_2(\om)+\frac{4}{\om^2} \al_{12} S_1(\om)\]
\[=\csc \om \bl(\frac{\al_{10}}{2} \sum_{r=0}^2\frac{b_{r,2}}{\om^r}-2\al_{11} \sum_{r=0}^1\frac{b_{r,1}}{\om^{r+1}}+\frac{4}{\om^2}\al_{12} b_{00}\br)=\csc \om \sum_{k=0}^2\frac{d_{1,k}}{\om^k},\]
where
\[d_{10}=\fs \al_{10}b_{02}=\fs(\cot^2 \om+\csc^2 \om),\qquad d_{11}=\fs \al_{10}b_{12}-2\al_{11}b_{01}=-\fs\cot \om,\]
\[d_{13}=\fs\al_{10}b_{22}-2\al_{11}b_{11}+4\al_{12}b_{00}=-\f{1}{6}.\]
In this manner we obtain after some effort
\[B_r(\om)=\csc \om \sum_{k=0}^{2r} \frac{d_{r,k}}{\om^k}\qquad (r\geq 1),\]
where
\[ d_{r,k}=\sum_{p=0}^k\frac{(-2)^{k-p}\al_{r,k-m}}{(p+2r-k)!} \,b_{m,p^*}\qquad (p^*=2r\!-\!k\!+\!p\!+\!1).\]

It then follows upon reversing the order of summation that
\[T_m(\om)=\csc \om \sum_{r=0}^{m-1} (\pi i/4)^r \sum_{k=0}^{2r} \frac{d_{r,k}}{\om^k}=\csc \om \sum_{k=0}^{2m-2} \frac{1}{\om^{k}} \sum_{r=\lceil k/2\rceil}^{m-1}(\pi i/4)^r d_{r,k}. \]
Finally, we obtain the correction term in the following form:
\begin{theorem}$\!\!\!.$\ \ The correction term $T_m(\om)$ appearing in
the expansion (\ref{e28}) has the form
\bee\label{e32}
T_m(\om)=\csc \om \sum_{k=0}^{2m-2} \frac{D_k(m)}{\om^k},
\ee
where the coefficients $D_k(m)\equiv D_k(m;\om)$ are given by
\bee\label{e33}
D_k(m)=\sum_{r=\lceil k/2\rceil}^{m-1}(\pi i/4)^r d_{r,k}.
\ee
\end{theorem}

The coefficients $D_k(\om)$ are easily seen to satisfy
\[D_k(m+1)=D_k(m)+(\pi i/4)^m d_{m,k}\qquad (0\leq k\leq 2m-4)\]
and involve the trigonometric functions $\csc \om$ and $\cot \om$.
We now present the explicit representation\footnote{The O$(\epsilon^3)$ contribution in $D_2(4)$ is zero.} of the coefficients $D_k(m)$ for $m=1, 2, 3$ and 4, where for brevity we have set $\epsilon=\pi i/4$.
\begin{eqnarray*}
m=1:\qquad&& D_0(1)=1;\\
\\
m=2:\qquad &&D_0(2)=D_0(1)+\fs \epsilon (\cot^2 \om+\csc^2 \om)\\
&&D_1(2)=-\fs \epsilon \cot \om,\ \ \ D_2(2)=-\f{1}{6}\epsilon;\\
\\
m=3:\qquad&& D_0(3)=D_0(2)+\f{1}{8}\epsilon^2(\cot^4 \om+18\cot^2 \om \csc^2 \om+5\csc^4 \om)\\
&&D_1(3)=D_1(2)-\f{5}{12}\epsilon^2 (\cot^3 \om+5\cot \om \csc^2 \om)\\
&&D_2(3)=D_2(2)-\f{5}{24}\epsilon^2 (\cot^2 \om+\csc^2 \om)\\
&&D_3(3)=-\f{1}{24} \epsilon^2 \cot \om,\ \ \ D_4(3)=\f{1}{72} \epsilon^2;\\
\\
m=4:\qquad&& D_0(4)=D_0(3)+\f{1}{48}\epsilon^3(\cot^6 \om+179 \cot^4 \om \csc^2 \om+479 \cot^2 \om \csc^4 \om+61 \csc^6 \om)\\
&&D_1(4)=D_1(3)-\f{7}{48}\epsilon^3(\cot^5 \om+58 \cot^3 \om \csc^2 \om+61\cot \om \csc^4 \om)\\
&&D_2(4)=D_2(3)\\
&&D_3(4)=D_3(3)+\f{7}{36}\epsilon^3 (\cot^3 \om +5\cot \om \csc^2 \om)\\
&&D_4(4)=D_4(3)+\f{49}{144}\epsilon^3 (\cot^2 \om+\csc^2 \om)\\
&&D_5(4)=\f{47}{144}\epsilon^3 \cot \om,\ \ \ D_6(4)=\f{1003}{6480} \epsilon^3.\\
\end{eqnarray*}
\vspace{0.6cm}

\begin{center}
{\bf 4.\ Numerical results}
\end{center}
\setcounter{section}{4}
\setcounter{equation}{0}
\renewcommand{\theequation}{\arabic{section}.\arabic{equation}}
In this section we describe computations using the expansion (\ref{e28}) with the correction term $T_m(\om)$ given in (\ref{e32}) and (\ref{e33}). The terms in the main sum in (\ref{e28}) decay rapidly beyond $n\simeq N_t$. In the neighbourhood of the transition point of $Q(\fs s,\pi n^2i)$ given by $\la_n=1$, we have
\[\eta_n\simeq \la_n-1=\frac{\pi^2}{\om^2}\bl(n^2-\frac{\om^2}{\pi^2}\br) \simeq \frac{2\pi i}{s}\bl(n^2-\frac{t}{2\pi}\br),\]
so that the argument of the modified complementary error function in (\ref{e28}) for large $t$ when $n\simeq N_t$ is
\[\fs\eta_n\sqrt{s}\simeq \pi t^{-1/2} e^{\pi i/4} (n^2-N_t^2)\simeq 2\pi t^{-1/2} e^{\pi i/4} N_t(n-N_t)\]
\[\simeq\sqrt{2\pi i}\,(n-N_t).\]
From the asymptotic behaviour
\[\mbox{erfc}\,(z;m)\sim\frac{e^{-z^2}}{\sqrt{\pi}}\,\frac{(-)^m(\fs)_m}{z^{2m+1}} \qquad(|z|\to\infty,\ \ |\arg\,z|<\f{3}{4}\pi),\]
it is then found that when $n=N_t+K$, where integer $K\ll N_t$, the magnitude of the terms in the main sum is given approximately by
\[|n^{-s}\mbox{erfc}\,(\fs\eta_n\sqrt{s};m)| =(4\pi K^2)^{-m-1/2} O((t/2\pi)^{-1/4}).\]
The decay of the terms is therefore controlled by $K^{-2m-1/2}$ together with a scaling factor depending weakly on $t$ like $t^{-1/4}$. Thus if $K=10$, for example, the magnitude of the terms in the main sum is approximately of order 
$10^{-14} t^{-1/4}$ when $m=4$.

\begin{table}[t]
\caption{\footnotesize{Details of the computations when $t=2600$, $N_t=20$, $Z(t)=-0.63210232$ for different values of $m$ and truncation index $K$ in the main sum.}}
\begin{center}
\begin{tabular}{|c|c|c|}
\hline
\mcol{3}{|c|}{$m=2$ \ \ Correction term=$+0.17012\ 33165\ 89694\ 33$} \\
\mcol{1}{|c|}{$K$} & \mcol{1}{|c|}{Main Sum} & \mcol{1}{|c|}{$|Z_{\mbox{\footnotesize{approx}}}-Z(t)|$}\\
[.1cm]\hline
&&\\[-0.25cm]
10 & $-0.46197\ 90063\ 18404\ 31$ & $2.269\times 10^{-08}$\\
20 & $-0.46197\ 90041\ 44898\ 25$ & $9.502\times 10^{-11}$\\
30 & $-0.46197\ 90040\ 66916\ 27$ & $1.704\times 10^{-11}$\\  
[.1cm]\hline

\mcol{3}{|c|}{$m=3$ \ \ Correction term=$+0.08577\ 29470\ 58489\ 99$} \\
\mcol{1}{|c|}{$K$} & \mcol{1}{|c|}{Main Sum} & \mcol{1}{|c|}{$|Z_{\mbox{\footnotesize{approx}}}-Z(t)|$}\\
[.1cm]\hline
&&\\[-0.25cm]
10 & $-0.54632\ 93736\ 04227\ 97$ & $2.315\times 10^{-11}$\\
20 & $-0.54632\ 93735\ 81347\ 73$ & $2.693\times 10^{-13}$\\
30 & $-0.54632\ 93735\ 81099\ 88$ & $2.146\times 10^{-14}$\\  
[.1cm]\hline

\mcol{3}{|c|}{$m=4$ \ \ Correction term=$-0.53258\ 74830\ 13204\ 32$} \\
\mcol{1}{|c|}{$K$} & \mcol{1}{|c|}{Main Sum} & \mcol{1}{|c|}{$|Z_{\mbox{\footnotesize{approx}}}-Z(t)|$}\\
[.1cm]\hline
&&\\[-0.25cm]
10 & $-1.16468\ 98036\ 52714\ 66$ & $5.792\times 10^{-14}$\\
20 & $-1.16468\ 98036\ 52772\ 40$ & $1.874\times 10^{-16}$\\
30 & $-1.16468\ 98036\ 52772\ 58$ & $6.799\times 10^{-18}$\\  
[.1cm]\hline
\end{tabular}
\end{center}
\end{table}
\begin{table}[t]
\caption{\footnotesize{Details of the computations when $t=2\times 10^{5}$, $N_t=178$, $Z(t)=-3.51142011$ for different values of $m$ and truncation index $K$ in the main sum.}}
\begin{center}
\begin{tabular}{|c|c|c|}
\hline
\mcol{3}{|c|}{$m=2$ \ \ Correction term=$+0.04418\ 05095\ 82215\ 20$} \\
\mcol{1}{|c|}{$K$} & \mcol{1}{|c|}{Main Sum} & \mcol{1}{|c|}{$|Z_{\mbox{\footnotesize{approx}}}-Z(t)|$}\\
[.1cm]\hline
&&\\[-0.25cm]
10 & $-3.46723\ 96042\ 56646\ 04$ & $5.980\times 10^{-10}$\\
20 & $-3.46723\ 96036\ 78939\ 00$ & $2.032\times 10^{-11}$\\
30 & $-3.46723\ 96036\ 61437\ 02$ & $2.815\times 10^{-12}$\\  
[.1cm]\hline

\mcol{3}{|c|}{$m=3$ \ \ Correction term=$+0.03215\ 07834\ 37290\ 07$} \\
\mcol{1}{|c|}{$K$} & \mcol{1}{|c|}{Main Sum} & \mcol{1}{|c|}{$|Z_{\mbox{\footnotesize{approx}}}-Z(t)|$}\\
[.1cm]\hline
&&\\[-0.25cm]
10 & $-3.47926\ 93298\ 09067\ 88$ & $5.521\times 10^{-12}$\\
20 & $-3.47926\ 93298\ 03596\ 26$ & $4.958\times 10^{-14}$\\
30 & $-3.47926\ 93298\ 03549\ 81$ & $3.130\times 10^{-15}$\\  
[.1cm]\hline

\mcol{3}{|c|}{$m=4$ \ \ Correction term=$-0.03047\ 11685\ 29774\ 29$} \\
\mcol{1}{|c|}{$K$} & \mcol{1}{|c|}{Main Sum} & \mcol{1}{|c|}{$|Z_{\mbox{\footnotesize{approx}}}-Z(t)|$}\\
[.1cm]\hline
&&\\[-0.25cm]
10 & $-3.54189\ 12817\ 70598\ 54$ & $1.251\times 10^{-14}$\\
20 & $-3.54189\ 12817\ 70611\ 01$ & $2.921\times 10^{-17}$\\
30 & $-3.54189\ 12817\ 70611\ 04$ & $8.377\times 10^{-19}$\\  
[.1cm]\hline
\end{tabular}
\end{center}
\end{table}

In the computations we define the truncated main sum as
\[\Re e^{i\vartheta(t)} \sum_{k=1}^{N_t+K} n^{-s} \mbox{erfc}\,(\fs\eta_n\sqrt{s};m)\]
\bee\label{e41}
=\Re e^{i\vartheta(t)} \bl\{\sum_{k=1}^{N_t} n^{-s} \bl(2-\mbox{erfc}\,(\fs\eta_n\sqrt{s};m)\br)+\sum_{k=N_t+1}^{N_t+K} n^{-s} \mbox{erfc}\,(\fs\eta_n\sqrt{s};m)\br\},
\ee
where\footnote{We note that the factor 2 in the sum over $1\leq n\leq N_t$ in (\ref{e41}) yields the standard Riemann-Siegel sum $2\sum_{n=1}^{N_t} n^{-1/2} \cos (\vartheta(t)-t\log\,n)$.} $N_t$ is the Riemann-Siegel cut-off value defined in (\ref{e12}) and we have made use of the result $\mbox{erfc}\,(-z;m)=2-\mbox{erfc}\,(z;m)$.
The correction term is given by
\[\Re \frac{e^{i\vartheta(t)}}{\sqrt{2i}}\,\bl(\frac{\pi e^{1/2}}{\om}\br)^{\!s} T_m(\om).\] 
The difference between these two contributions then yields the value $Z_{\mbox{\footnotesize{approx}}}$.
An example of the results is displayed in Tables 1 and 2, where the value $Z(t)$ was obtained by {\it Mathematica} using the command {\tt RiemannSiegelZ[t]}..
The values in the tables confirm that the accuracy increases as both $m$ and the truncation index $K$ increase.
\vspace{0.6cm}

\begin{center}
{\bf 5.\ Concluding remarks}
\end{center}
\setcounter{section}{5}
\setcounter{equation}{0}
\renewcommand{\theequation}{\arabic{section}.\arabic{equation}}
We have revisited an expansion derived in \cite{PC} for $Z(t)$ on the critical line $s=\fs+it$ as $t\to+\infty$ in
which the main sum is the original Dirichlet series smoothed by a complementary error function. The correction term in this expansion has been expressed as a series in inverse powers of $\om$, where $\om^2=\pi s/(2i)$, multiplied by coefficients involving $\csc \om$ and its derivatives. Numerical results are presented to illustrate the accuracy achievable with this expansion.

A difficulty arises in the use of the expansion (\ref{e28}) when $(t/2\pi)^{1/2}$ is close to an integer i.e., at a discontinuity in $N_t$. This arises because
\[\frac{\om}{\pi}=\bl(\frac{t}{2\pi}\br)^{\!1/2} \bl(1-\frac{i}{2t}\br)^{\!1/2},\]
so that $\csc \om$ and $\cot \om$ become large, but never singular as $\om$ always has a small imaginary part (that decreases with increasing $t$). This results in the term in the main sum and correction term corresponding to $\eta_n\simeq 0$ becoming large. A means of overcoming this problem, by removing and combining these terms, is given in \cite[Section 5]{PC}.

\end{document}